\documentclass[oneside,reqno]{amsart}
\usepackage{amsmath}
\usepackage[mathscr]{euscript}
\newtheorem{theorem}{Theorem}[section]
\newtheorem{lemma}[theorem]{Lemma}

\theoremstyle{definition}
\newtheorem{definition}[theorem]{Definition}

\numberwithin{equation}{section}

\newcommand{\blankbox}[2]

\begin{document}
\title{Some results for the commutators of Generalized Hausdorff operator}
\author{Amjad Hussain$^{1,*}$}
\author{Amna Ajaib$^1$}
%\date{February, 2013}
\footnote{$^{}$Department of Mathematics, Quaid-I-Azam University 45320, Islamabad 44000, Pakistan\\
$^*$Corresponding Author: ahabbasi123@yahoo.com (A. Hussain)}


\begin{abstract}
In this paper, we study commutator of generalized Hausdorff operator on function spaces. We mainly discuss the continuity criteria for such commutator operator when the symbol functions are either from central-$BMO$ or Lipschitz class of functions.

\end{abstract}

\maketitle

%\footnote{The second author is the corresponding author.}%

\section{\textbf{Introduction}}

 To study the continuity properties of Hausdorff operators on function spaces, an appropriate point for opening discussion is the one dimensional Hausdorff operator:
\begin{eqnarray}\label{E1}\mathscr{H}_\Phi f(x)=\int_0^\infty\frac{\Phi(t)}{t}f(\frac{x}{t})dt,\qquad x\in \mathbb R,\end{eqnarray}
where $\Phi\in L^1( \mathbb R).$ Also, in the interest of convenience, it is usually assumed that functions $f$ are initially in Schwartz class $\mathcal S(\mathbb R).$ In fact, there are many popular and important operators in analysis which become special cases of $\mathscr{H}_\Phi,$ if the function $\Phi$ is suitably chosen. Among many others, here we mention a few namely, the Ces\`{a}ro operator \cite{K}, the Hardy operator and its adjoint operator \cite{CG}.
The survey articles \cite{CFW,EL1} may serve as an important source of information about history and recent advances in the study of this operator.

In the case of $n$-dimensional Euclidian space $\mathbb R^n(n\ge2),$ an extension of (\ref{E1}) was studied in \cite{LL} and is given by
\begin{eqnarray}\label{E2}H_{\Phi,A}f(x)=\int_{\mathbb R^n} \frac{\Phi(y)}{|y|^n} f( A(y)x) dy,\end{eqnarray}
where $A(y)$ is an $n$-th order square matrix, which satisfy $\det A(y)\ne0$ almost everywhere in the support of $\Phi\in L_{\rm loc}^1(\mathbb R^n).$ Using duality approach, Lerner and Liflayand \cite{LL} obtained the boundedness of $H_{\Phi, A}$ on real Hardy space $H^1(\mathbb R^n),$ after which the same problem was reconsidered in \cite{CZ,EL2,RF2}.

In the above definition (\ref{E2}), if one considers $A(y)=\text{diag}[1/|y|,1/|y|,...,1/|y|],$ then we obtain a natural $n$-dimensional extension of (\ref{E1}) given by
\begin{eqnarray}\label{E3}H_{\Phi}f(x)=\int_{\mathbb R^n} \frac{\Phi(y)}{|y|^n} f(\frac{x}{|y|}) dy.\end{eqnarray}

Like one dimensional Hausdorff operator \cite{LM,LM1}, the operators $H_{\Phi}$ and $H_{\Phi,A}$ also enjoy a variety of estimates on different function spaces; see \cite{BL,CFL,CHZ,GWAG,GZ,GZ1,KP,HKQ,HA,HG2,RF1,RFW} and the references therein.

Likewise the study of boundedness properties of commutators operators on function spaces is another most important problem in analysis. The study can be utilized for the characterization of function spaces, in the regularity theory to special kind of partial differential equations (PDEs) and in the well posedness problems of solution to many PDEs. Therefore, the study of boundedness results for the commutators of Hausdorff operators is as important as the study of Hausdorff operators itself. In an exploratory research, one can find a very few papers discussing the boundedness of commutators of various Hausdorff operators \cite{GJ,HA1,HG1,HG2,XW}, except that of $H_{\Phi,A}.$ Recently, in \cite{HA}, the author and his collaborator defined the commutator generated by generalized Hausdorff operator $H_{\Phi,A}$ and locally integrable function $b$ as:
\begin{eqnarray*}H_{\Phi,A}^b(f)(x)=b(x)H_{\Phi,A}(f)(x)-H_{\Phi,A}(bf)(x),\end{eqnarray*}
and constructed weighted estimates for it on central Morrey spaces, when $b\in C\dot{M}O(\mathbb R^n).$ Also, they raised an open question regarding $L^p(\mathbb R^n)$ boundedness of $H_{\Phi,A}^b$.
Here, we give partially positive answer to this question by establishing $L^p(\mathbb R^n)\rightarrow L^q(\mathbb R^n)$ estimates when the symbol function belongs to the Lipschitz class of functions. However, in the case $b\in C\dot{M}O(\mathbb R^n),$ the question of $L^p(\mathbb R^n)$ boundedness of $H_{\Phi,A}^b$ still remains open.

In this paper, we give Lipschitz estimates for the commutators of $H_{\Phi,A}$ on Lebesgue, Morrey and Herz-Type spaces and thus generalize some results presented in \cite{GZ,HA1}. In addition, when $b\in\dot{\Lambda}_\beta(\mathbb R^n),$ we obtain $L^p\rightarrow \dot{F}^{\beta,\infty}_p$ boundedness for $H_{\Phi,A}^b,$ where $\dot{F}^{\beta,\infty}_p$ is the homogeneous Triebel-Lizorkin space. Also, we estimate $H_{\Phi,A}^b$ on Herz-type spaces when $b\in C\dot{M}O(\mathbb R^n).$ The significance of our results lies in the fact that Herz-Type spaces are used in the characterization of multiplier on the Hardy spaces \cite{BS} and in the study of certain kind of PDEs \cite{LY}.

 The plot of this paper is as follows. The second section contains some basic definitions and notations likewise some necessary lemmas which will be used in the succeeding sections. Lipschitz estimates for $H_{\Phi,A}^b$ are established in the third section. While results regarding central-$BMO$ estimates for $H_{\Phi,A}^b$ on Herz-type spaces are stated and proved in the last section.

\section{\textbf{Preliminaries}}

In 1938, Morrey \cite{CM} studied the local behavior of second order elliptic and parabolic PDEs and introduced a function space what is called Morrey space.
\begin{definition} Let $1\leq p<\infty,~~~~~0\leq\lambda\leq  n.$ The Morrey space $L^{p,\lambda}(\mathbb R^n)$ is the set of all functions $f$
$$L^{p,\lambda}(\mathbb R^n)=\left\{f \in L_{\rm loc} ^{p}(\mathbb{R}^{n}):\|f\|_{L^{p,\lambda}(\mathbb R^n)}< \infty\right\},$$
satisfying $$\|f\|_{L^{p,\lambda}(\mathbb R^n)}=\sup_{r >0,x_0 \in \Bbb{R}^{n}}\left(\frac{1}{r^{\lambda}}
\int_{Q(x_0,r)}|f(x)|^{p}dx\right)^{1/p},$$
where $Q=Q(x_0,r)$ denotes the cube centered at $x_0$ with side length $r$ along the coordinate axes.
\end{definition}
It is easy to see that $L^{p,0}(\mathbb R^n)=L^{p}(\mathbb R^n)$
and $L^{p,n}(\mathbb R^n)=L^{\infty}(\mathbb R^n).$ If $n<\lambda,$ then we have
$L^{p,\lambda}(\mathbb R^n)={0}$. Hence, we only consider the case $0<\lambda<n.$

Let $B_k:=\{x\in \mathbb R^n:|x|<2^k\}, C_k=B_k/B_{k-1}$ for $k\in \mathbb Z.$ Then we will consider the following definition of homogeneous Herz space.
\begin{definition}
Let $\alpha\in\mathbb R, 0<p,q<\infty.$ The homogeneous Herz space $\dot{K}^{\alpha,p}_q(\mathbb R^n)$ is defined by
$$\dot{K}^{\alpha,p}_q(\mathbb R^n):=\left\{f \in L_{\rm loc}^q(\mathbb H^{n}/\{0\}):\|f\|_{\dot{K}^{\alpha,p}_q(\mathbb R^n)}< \infty\right\},$$
where
$$\|f\|_{\dot{K}^{\alpha,p}_q(\mathbb R^n)}=\left\{\sum_{k=-\infty}^\infty 2^{k\alpha p}\|f\|_{L^q(C_k)}^p\right\}^{1/p}.$$
\end{definition}
 It is easy to verify that $\dot{K}^{0,p}_p(\mathbb R^n)=L^p(\mathbb R^n).$ Hence, Herz space can be considered as an extension of Lebesgue space $L^p(\mathbb R^n).$

Similarly the homogenous Herz-Morrey space $M\dot{K}^{\alpha,\lambda}_{p,q}(\mathbb R^n)$ is given by the following definition.
\begin{definition}\label{D2}
Let $\alpha\in\mathbb R, 0<p,q<\infty, \lambda\ge0.$ The homogeneous Herz-Morrey space $M\dot{K}^{\alpha,\lambda}_{p,q}(\mathbb R^n)$ is defined by
$$M\dot{K}^{\alpha,\lambda}_{p,q}(\mathbb R^n):=\left\{f \in L_{\rm loc}^q(\mathbb R^{n}/\{0\}):\|f\|_{M\dot{K}^{\alpha,\lambda}_{p,q}(\mathbb R^n)}< \infty\right\},$$
where
$$\|f\|_{M\dot{K}^{\alpha,\lambda}_{p,q}(\mathbb R^n)}=\sup_{k_0\in\mathbb Z}2^{-k_0\lambda}\left\{\sum_{k=-\infty}^{k_0} 2^{k\alpha p}\|f\|_{L^q(C_k)}^p\right\}^{1/p}.$$
\end{definition}
Obviously, $M\dot{K}^{\alpha,0}_{p,q}(\mathbb R^n)=\dot{K}^{\alpha,p}_{q}(\mathbb R^n)$ and
$L^{q,\lambda}(\mathbb R^n)\subset M \dot{K}_{q,q}^{0,\lambda}(\mathbb R^n).$

\begin{definition}[] Let $1<q<\infty$. A function
$f\in L_{loc}^q(\mathbb R^n)$ is said to belong to the space ${C\dot{M}O}^{q}(\mathbb{R}^n)$ if
$$\|f\|_{{C\dot{M}O}^{q}(\mathbb{R}^{n})}=
\sup\limits_{r>0}\left(\frac{1}{|B(0,r)|}\int_{B(0,r)}|f(x)-f_{B(0,r)}|^qdx\right)^{1/q}<\infty,$$
where $f_B=\frac{1}{|B|}\int_Bf(x)dx$ denote the average value of $f$ over $B$. Obviously, we have $BMO(\mathbb{R}^n)\subset C\dot{M}O^{q}(\mathbb{R}^n) $
for all $ 1 \leq q < \infty ,$ and $ {\rm
C\dot{M}O}^{q}(\mathbb{R}^n) \subset {\rm
C\dot{M}O}^{p}(\mathbb{R}^n) $, $ 1 \leq p < q< \infty $.
\end{definition}
\begin{definition}[] Let $0<\beta<1.$ The
Lipschitz space $\dot{\Lambda}_{\beta}(\mathbb{R}^{n})$ is
defined by
$$\|f\|_{\dot{\Lambda}_{\beta}(\mathbb{R}^{n})}=\sup_{x,h\in\mathbb{R}^{n}}\frac{|f(x+h)-f(x)|}{|h|^{\beta}}<\infty.$$
\end{definition}

Next,  for a non-singular matrix $B,$ we consider the following definition of norm
\begin{eqnarray}\label{E21}
\|B\|=\sup_{x\in \mathbb R^n,x\neq0}\frac{|Bx|}{|x|},\end{eqnarray}
which implies that
\begin{eqnarray}\label{E22}\begin{split}
\|B\|^{-n}\le|\det(B^{-1})|\le\|B^{-1}\|^{n}.\end{split}\end{eqnarray}

Finally, for $0\leq\beta<n,$ we
define the fractional maximal function as
$$M_\beta f(x)=\sup\limits_{Q\ni x}\frac{1}{|Q|^{1-\frac{\beta}{n}}}\int_Q|f(y)|dy,$$
where the supremum is taken over all cubes $Q$ containing $x$.
Notice that when $\beta=0,$ we obtain the usual Hardy Littlewood
maximal function $M=M_0.$

This finishes the streak of definitions concerning function spaces, norm of a matrix and maximal function.
We are now in position to state some lemmas which will be used to prove our main results.
\begin{lemma}\label{LL}(\cite{PM}) Let $0<\beta<1$ and $f\in\dot{\Lambda}_\beta(\mathbb R^n),$ then
 for any cube $Q\subset\mathbb R^n,$
$$\sup_{x\in Q}|f(x)-f_Q|\leq C|Q|^{\frac{\beta}{n}}\|f\|_{\dot{\Lambda}_\beta(\mathbb R^n)}.$$
\end{lemma}

\begin{lemma}\label{LA} Let $0<\beta<1$ and $b\in\dot{\Lambda}_\beta(\mathbb R^n),$ then
\begin{eqnarray*}\begin{split} &M(H_{\Phi,A}^bf)(x)\\& \le C\|b\|_{\dot{\Lambda}_\beta(\mathbb R^n)}\int_{\mathbb R^n}\frac{|\Phi(y)|}{|y|^n}\max\{1,|\det A^{-1}(y)|^{\beta/n}\}(1+\|A(y)\|^\beta)M_\beta(f)(A(y)x)dy.\end{split}\end{eqnarray*}
\end{lemma}
\noindent Proof: Let us consider a cube $Q\subset\mathbb R^n$ such that $x\in Q$ and
\begin{eqnarray*}    \begin{aligned}[b]
\frac{1}{|Q|}\int_{Q}|H_{\Phi,A}^bf(z)|dz=&\frac{1}{|Q|}\int_{Q}\left|\int_{\mathbb R^n}\frac{\Phi(y)}{|y|^n}(b(z)-b(A(y) z))f(A(y) z)dy\right|dz\\
&\le\frac{1}{|Q|}\int_{\mathbb R^n}\frac{|\Phi(y)|}{|y|^n}\int_{Q}\left|(b(z)-b(A(y) z))f(A(y) z)\right|dzdy\\
&\le\frac{1}{|Q|}\int_{\mathbb R^n}\frac{|\Phi(y)|}{|y|^n}\int_{Q}\left|(b(z)-b_{Q})f(A(y) z)\right|dzdy\\
&+\frac{1}{|Q|}\int_{\mathbb R^n}\frac{|\Phi(y)|}{|y|^n}\int_{Q}\left|(b_{Q}-b_{A(y)Q})f(A(y) z)\right|dzdy\\
&+\frac{1}{|Q|}\int_{\mathbb R^n}\frac{|\Phi(y)|}{|y|^n}\int_{Q}\left|(b(A(y) z)-b_{A(y)Q})f(A(y) z)\right|dzdy\\
&=I_1+I_2+I_3.
\end{aligned}\end{eqnarray*}

In the approximation of $I_1,$ we use Lemma \ref{LL} to obtain
\begin{eqnarray*}\begin{split}I_1&\leq C\|b\|_{\dot{\Lambda}_\beta(\mathbb R^n)}\int_{\mathbb R^n}\frac{|\Phi(y)|}{|y|^n}
\left(\frac{1}{|Q|^{1-\beta/n}}\int_Q |f(A(y)z)| dz\right)dy\\&\leq C\|b\|_{\dot{\Lambda}_\beta(\mathbb R^n)}\int_{\mathbb R^n}\frac{|\Phi(y)|}{|y|^n}
\left(\frac{1}{|A(y)Q|^{1-\beta/n}}\int_{A(y)Q} |f(z)| dz\right)|\det A^{-1}(y)|^{\beta/n}dy\\
&\leq C\|b\|_{\dot{\Lambda}_\beta(\mathbb R^n)}\int_{\mathbb R^n}\frac{|\Phi(y)|}{|y|^n}|\det A^{-1}(y)|^{\beta/n}
 M_{\beta}( f)(A(y)x)dy.\end{split}\end{eqnarray*}

In order to estimate $I_2,$ we have to approximate $|b_Q-b_{A(y)Q}|.$ For $0<\beta<1,$ we obtain
\begin{eqnarray*}\begin{split}|b_Q-b_{A(y)Q}|&\leq\frac{1}{|Q|}\int_{Q}|b(z)-b_{A(y)Q}|dz\\
&\leq\frac{1}{|Q|}\frac{1}{|A(y)Q|}\int_Q\int_{A(y)Q}|b(z)-b(t)|dtdz\\
&\leq\|b\|_{\dot{\Lambda}_\beta(\mathbb R^n)}
\left(\frac{1}{|Q|}\int_Q|x|^{\beta}dz+\frac{1}{|A(y)Q|}\int_{A(y)Q}|t|^{\beta}dt\right)\\
&\leq C|Q|^{\beta/n}\|b\|_{\dot{\Lambda}_\beta(\mathbb R^n)}(1+\|A(y)\|^{\beta})\end{split}
\end{eqnarray*}
Therefore.
\begin{eqnarray*}\begin{split}I_2&\leq C\|b\|_{\dot{\Lambda}_\beta(\mathbb R^n)}
\int_{\mathbb R^n}\frac{|\Phi(y)|}{|y|^n}(1+\|A(y)\|^{\beta})
\left(\frac{1}{|Q|^{1-\beta/n}}\int_Q |f(x A(y))| dz\right)dy\\
&\leq C\|b\|_{\dot{\Lambda}_\beta(\mathbb R^n)}\int_{\mathbb R^n}\frac{|\Phi(y)|}{|y|^n}(1+\|A(y)\|^{\beta})
|\det A^{-1}(y)|^{\beta/n}
 M_{\beta}( f)(A(y)x)dy.\end{split}\end{eqnarray*}

 It remains to estimate $I_3.$ By virtue of Lemma \ref{LL}, we get
 \begin{eqnarray*}\begin{split}I_3&=\frac{1}{|Q|}\int_{\mathbb R^n}\frac{|\Phi(y)|}{|y|^n}\int_{Q}\left|(b(A(y) z)-b_{A(y)Q})f(A(y) z)\right|dzdy\\
 &=\frac{1}{|Q|}\int_{\mathbb R^n}\frac{|\Phi(y)|}{|y|^n}|\det A(y)|^{-1}\int_{A(y)Q}|b( z)-b_{A(y)Q}||f(z)|dzdy\\
 &\leq C\|b\|_{\dot{\Lambda}_\beta(\mathbb R^n)}\int_{\mathbb R^n}\frac{|\Phi(y)|}{|y|^n}
\left(\frac{1}{|A(y)Q|^{1-\beta/n}}\int_{A(y)Q} |f(z)| dz\right)dy\\
&\leq C\|b\|_{\dot{\Lambda}_\beta(\mathbb R^n)}\int_{\mathbb R^n}\frac{|\Phi(y)|}{|y|^n}
 M_{\beta}( f)(A(y)x)dy.\end{split}\end{eqnarray*}

We combine the estimates for $I_1,$ $I_2$ and $I_3,$ to have
\begin{eqnarray*}\begin{split} &\frac{1}{|Q|}\int_{Q}|H_{\Phi,A}^bf(z)|dz\\& \le C\|b\|_{\dot{\Lambda}(\mathbb R^n)}\int_{\mathbb R^n}\frac{\Phi(y)}{|y|^n}\max\{1,|\det A^{-1}(y)|^{\beta/n}\}(1+\|A(y)\|^\beta)M_\beta(f)(A(y)x).\end{split}\end{eqnarray*}
Finally, taking supremum over all $Q$ such that $x\in Q,$ we get the desired result.
\begin{lemma}(\cite{BGG})\label{LM}
Let $0<\beta<n,$ $1<p<n/\beta,$ $0<\lambda<n-\beta p$ and $1/q=1/p-\beta/(n-\lambda)$. Then $M_\beta$
is bounded from $L^{p,\lambda}({\mathbb R^n})$ to $L^{q,\lambda}({\mathbb R^n}).$
\end{lemma}
\begin{lemma}(\cite{LYH})\label{LH} Let $0<p_1\le\infty,p_1\le p_2\le\infty, 0<\beta<1<q_1<n/\beta$ and $1/q_1-1/q_2=\beta/n.$ If $-n/q_1+\beta<\alpha<n(1-1/q_1),$ then
$$\|M_\beta f\|_{\dot{K}^{\alpha,p_2}_{q_2}(\mathbb R^n)}\le C\|f\|_{\dot{K}^{\alpha,p_1}_{q_1}(\mathbb R^n)}.$$\end{lemma}
\begin{lemma}(\cite{GU})\label{LMH} Let $0<p_1\le\infty,p_1\le p_2\le\infty, 0<\beta<1<q_1<n/\beta,1/q_1-1/q_2=\beta/n$ and $\lambda>0.$ If $-n/q_1+\beta+\lambda<\alpha<n(1-1/q_1)+\lambda,$ then
$$\|M_\beta f\|_{M\dot{K}^{\alpha,\lambda}_{p_2,q_2}(\mathbb R^n)}\le C\|f\|_{M\dot{K}^{\alpha,\lambda}_{p_1,q_1}(\mathbb R^n)}.$$\end{lemma}
\begin{lemma}(\cite{PM})\label{LT} Let $0<\beta<1<p<\infty,$  then
$$\|f\|_{F^{\beta,\infty}_p(\mathbb R^n)}\approx\left\|\sup\limits_{Q\ni x}\frac{1}{|Q|^{1+\beta/n}}\int_{Q}|f(z)-f_Q|dz\right\|_{L^p(\mathbb R^n)}.$$\end{lemma}

\section{\textbf{Lipschitz Estimates for $H_{\Phi,A}^b$ on Function Spaces}}
As we stated in the introduction, this section is centered on obtaining estimates for $H_{\Phi,A}^b$ on function spaces. In this regard, our main results are as below.
\subsection{Main Results}
\begin{theorem}\label{T2}
Let $0<\beta<1<p<n/\beta,$ $0<\lambda<n-\beta p$ and $1/q=1/p-\beta/(n-\lambda)$. If $b\in \dot{\Lambda}_\beta(\mathbb R^n),$ then
\begin{eqnarray*}\begin{split}&\|H_{\Phi,A}^bf\|_{L^{q,\lambda}({\mathbb R^n})}\le C K_1 \|b\|_{\dot{\Lambda}_\beta(\mathbb R^n)}\|f\|_{L^{p,\lambda}({\mathbb R^n})} ,\end{split}\end{eqnarray*}
where $K_1$ is
\begin{eqnarray*}\begin{split} \int_{\mathbb R^n}\frac{|\Phi(y)|}{|y|^n}\max\{|\det A^{-1}(y)|^{1/q-\lambda/q},|\det A^{-1}(y)|^{\beta/n+1/q-\lambda/q}\}(1+\|A(y)\|^\beta)dy.\end{split}\end{eqnarray*}
\end{theorem}
\begin{theorem}\label{T1}
Let $0<\beta<1<p<n/\beta,$ and $1/q=1/p-\beta/n$. If $b\in \dot{\Lambda}_\beta(\mathbb R^n),$ then
\begin{eqnarray*}\begin{split}&\|H_{\Phi,A}^bf\|_{L^{q}({\mathbb R^n})}\le C K_2\|b\|_{\dot{\Lambda}_\beta(\mathbb R^n)}\|f\|_{L^{p}({\mathbb R^n})},\end{split}\end{eqnarray*}
where $K_1$ is
\begin{eqnarray*}\begin{split}\int_{\mathbb R^n}\frac{|\Phi(y)|}{|y|^n}\max\{|\det A^{-1}(y)|^{1/q},|\det A^{-1}(y)|^{1/p}\}(1+\|A(y)\|^\beta)dy,\end{split}\end{eqnarray*}
\end{theorem}
\begin{theorem}\label{T4}
Let $0<p_1\le\infty,p_1\le p_2\le\infty, 0<\beta<1<q_1<n/\beta,1/q_1-1/q_2=\beta/n$ and $\lambda>0.$ If $-n/q_1+\beta+\lambda<\alpha<n(1-1/q_1)+\lambda$ and $b\in \dot{\Lambda}_\beta(\mathbb R^n),$ then
\begin{eqnarray*}\begin{split}&\|H_{\Phi,A}^bf\|_{M\dot{K}^{\alpha,\lambda}_{p_2,q_2}(\mathbb R^n)}\le C K_3\|b\|_{\dot{\Lambda}_\beta(\mathbb R^n)}\|f\|_{M\dot{K}^{\alpha,\lambda}_{p_1,q_1}(\mathbb R^n)},\end{split}\end{eqnarray*}
where $K_3$ is
\begin{eqnarray*}\begin{split}\int_{\mathbb R^n}\frac{|\Phi(y)|}{|y|^n}\max\{|\det A^{-1}(y)|^{1/q_2},|\det A^{-1}(y)|^{1/q_1}\}(1+\|A(y)\|^\beta)G_{\alpha,\lambda}(y)dy,\end{split}\end{eqnarray*}
and
\begin{eqnarray}\label{GLE}G_{\alpha,\lambda}(y)= \begin{cases}1+\log_2(\|A(y)\|\|A^{-1}(y)\|),  & \alpha=\lambda, \\
\|A^{-1}(y)\|^{\alpha-\lambda},  & \alpha>\lambda, \\
\|A(y)\|^{\lambda-\alpha},  & \alpha<\lambda.
\end{cases}\end{eqnarray}
\end{theorem}
\begin{theorem}\label{T3}
Let $0<p_1\le\infty,p_1\le p_2\le\infty, 0<\beta<1<q_1<n/\beta$ and $1/q_1-1/q_2=\beta/n.$ If $-n/q_1+\beta<\alpha<n(1-1/q_1)$ and $b\in \dot{\Lambda}_\beta(\mathbb R^n),$ then
\begin{eqnarray*}\begin{split}&\|H_{\Phi,A}^bf\|_{\dot{K}^{\alpha,p_2}_{q_2}(\mathbb R^n)}\le C K_4\|b\|_{\dot{\Lambda}_\beta(\mathbb R^n)}\|f\|_{\dot{K}^{\alpha,p_1}_{q_1}(\mathbb R^n)},\end{split}\end{eqnarray*}
where $K_4$ is
\begin{eqnarray*}\begin{split}\int_{\mathbb R^n}\frac{|\Phi(y)|}{|y|^n}\max\{|\det A^{-1}(y)|^{1/q_2},|\det A^{-1}(y)|^{1/q_1}\}(1+\|A(y)\|^\beta)\widetilde{G}_{\alpha}(y)dy,\end{split}\end{eqnarray*}
and
\begin{eqnarray}\label{GE} \widetilde{G}_{\alpha}(y)=\begin{cases}1+\log_2(\|A(y)\|\|A^{-1}(y)\|),  & \alpha=0, \\
\|A^{-1}(y)\|^{\alpha},  & \alpha>0, \\
\|A(y)\|^{-\alpha},  & \alpha<0.
\end{cases}\end{eqnarray}
\end{theorem}
\begin{theorem}\label{T5} Let $0<\beta<1<p<\infty$ and $b\in \dot{\Lambda}_\beta(\mathbb R^n),$ then
\begin{eqnarray*}\begin{split}\|H_{\Phi,A}^bf\|_{F^{\beta,\infty}_p(\mathbb R^n)}\le C K_5\|b\|_{\dot{\Lambda}_\beta(\mathbb R^n)}\|f\|_{L^{p}({\mathbb R^n})},\end{split}\end{eqnarray*}
where $K_5$ is
\begin{eqnarray*}\begin{split}\int_{\mathbb R^n}\frac{|\Phi(y)|}{|y|^n}\max\{|\det A^{-1}(y)|^{-1/p},|\det A^{-1}(y)|^{1+\beta/n-1/p}\}(1+\|A(y)\|^\beta)dy.\end{split}\end{eqnarray*}
\end{theorem}

\subsection{Proof of the Main Results}:\\

 Since Lemma \ref{LA} will be used frequently in the proofs of our main results for this section, therefore, we sometimes use the following notation $$\phi(y)=\frac{|\Phi(y)|}{|y|^n}\max\{1,|\det A^{-1}(y)|^{\beta/n}\}(1+\|A(y)\|^\beta),$$
for our convenience.\\

\noindent\textbf{\textit{Proof of Theorem \ref{T2}}}: In view of Lemma \ref{LA} and the Minkowski inequality, we have
\begin{eqnarray*}\begin{split} \|M(H_{\Phi,A}^bf)(\cdot)\|_{L^{q,\lambda}(\mathbb R^n)} \le C\|b\|_{\dot{\Lambda}_\beta(\mathbb R^n)}\int_{\mathbb R^n}\phi(y)\|M_\beta(f)(A(y)\cdot)\|_{L^{q,\lambda}(\mathbb R^n)}dy.\end{split}\end{eqnarray*}

Using the scaling argument and the fact that $|H_{\Phi,A}^bf(x)|\le M(H_{\Phi,A}^bf)(x)$ $a.e.,$ we get
\begin{eqnarray*}\begin{split} \|H_{\Phi,A}^bf\|_{L^{q,\lambda}(\mathbb R^n)}
 \le C\|b\|_{\dot{\Lambda}_\beta(\mathbb R^n)}\|M_\beta f\|_{L^{p,\lambda}(\mathbb R^n)}\int_{\mathbb R^n}\phi(y)|\det A^{-1}(y)|^{1/q-\lambda/q}dy.\end{split}\end{eqnarray*}
 Lastly, we apply Lemma \ref{LM} to obtain the desired result.\\

\noindent\textbf{\textit{Proof of Theorem \ref{T1}:}} Following the same procedure as followed in the proof of Theorem \ref{T2}, the proof of this Theorem can be easily obtained.\\

\noindent\textbf{\textit{Proof of Theorem \ref{T4}:}}
 Making use of Lemma \ref{LA} and the Minkowski inequality, one has
\begin{eqnarray}\label{E36}\begin{aligned}[b] &\|M(H_{\Phi,A}^bf)(\cdot)\|_{M\dot{K}^{\alpha,\lambda}_{p_2,q_2}(\mathbb R^n)}&\\
&\le C\|b\|_{\dot{\Lambda}_\beta(\mathbb R^n)}
\int_{\mathbb R^n}\phi(y)\|M_\beta f(A(y)\cdot)\|_{M\dot{K}^{\alpha,\lambda}_{p_2,q_2}(\mathbb R^n)}dy\\
&=C\|b\|_{\dot{\Lambda}_\beta(\mathbb R^n)}
\int_{\mathbb R^n}\phi(y)\sup_{k_{0}\in \mathbb Z}2^{-k_0\lambda}
\left\{\sum^{k_{0}}_{k=-\infty} 2^{k\alpha p}\|M_\beta f(A(y)\cdot) \|_{L^{q}(C_k)}^p\right\}^{1/p}dy\\
&=C\|b\|_{\dot{\Lambda}_\beta(\mathbb R^n)}
\int_{\mathbb R^n}\phi(y)|\det A^{-1}(y)|^{1/{q_2}}\\
&\qquad\times\sup_{k_{0}\in \mathbb Z}2^{-k_0\lambda}
\left\{\sum^{k_{0}}_{k=-\infty} 2^{k\alpha p}\|M_\beta f(\cdot) \|_{L^{q}(A(y)C_k)}^p\right\}^{1/p}dy.
\end{aligned}\end{eqnarray}

In order to estimate $\|M_\beta f(\cdot) \|_{L^{q}(A(y)C_k)},$ we follow the method used in \cite{RFW}.
Thus, by definition of $C_k$ and (\ref{E22}) it is easy to see that
\begin{eqnarray*}A(y)C_k\subset\{x:\|A^{-1}(y)\|^{-1}2^{k-1}<|x|<\|A(y)\|2^k\}.\end{eqnarray*}

Next, for any $y$ in the support of $ \Phi$ there exist an integer $l$  such that
\begin{eqnarray}\label{E31}2^{l}<\|A^{-1}(y)\|^{-1}<2^{l+1}.\end{eqnarray}
Furthermore, the relation $\|A^{-1}(y)\|^{-1}\le\|A(y)\|$ implies the existence of non-negative integer $m$ such that
\begin{eqnarray}\label{E32}2^{l+m}<\|A(y)\|<2^{l+m+1}.\end{eqnarray}

Inequalities (\ref{E31}) and (\ref{E32}) define the bounds for $m,$ that is
\begin{eqnarray*}\log_2(\|A(y)\|\|A^{-1}(y)\|/2)<m<\log_2(2\|A(y\|\|A^{-1}(y))\|),\end{eqnarray*}
and lead us to have
\begin{eqnarray*}A(y)C_k\subset\{x:2^{l+k-1}<|x|<2^{k+l+m+1}\}.\end{eqnarray*}
Hence,
\begin{eqnarray}\label{EA}\begin{aligned}
\|M_\beta f(\cdot)\|_{L^{q_2}(A(y)C_k)}\le\sum_{j=l}^{l+m+1}\|M_\beta f(\cdot)\|_{L^{q_2}(C_{k+j})}.
\end{aligned}
\end{eqnarray}

Incorporating the inequality (\ref{EA}) into (\ref{E36}), we obtain
\begin{eqnarray}\begin{aligned}[b]\label{E33}& \|M(H_{\Phi,A}^bf)(\cdot)\|_{M\dot{K}^{\alpha,\lambda}_{p_2,q_2}(\mathbb R^n)}&\\
&\le C\|b\|_{\dot{\Lambda}_\beta(\mathbb R^n)}
\int_{\mathbb R^n}\phi(y)|\det A^{-1}(y)|^{1/{q_2}}\\
&\quad\times\sup_{k_{0}\in \mathbb Z}2^{-k_0\lambda}
\left\{\sum^{k_{0}}_{k=-\infty} \left(\sum_{j=l}^{l+m+1}2^{k\alpha}\|M_\beta f(\cdot)\|_{L^{q_2}(C_{k+j})}\right)^p\right\}^{1/p}dy\\
&\le C\|b\|_{\dot{\Lambda}_\beta(\mathbb R^n)}
\int_{\mathbb R^n}\phi(y)|\det A^{-1}(y)|^{1/{q_2}}\\
&\quad\times\sup_{k_{0}\in \mathbb Z}2^{-k_0\lambda}
\sum_{j=l}^{l+m+1}2^{-j\alpha}\left\{\sum^{k_{0}+j}_{k=-\infty} 2^{k\alpha p }\|M_\beta f(\cdot)\|_{L^{q_2}(C_{k})}^p\right\}^{1/p}dy\\
&\le C\|b\|_{\dot{\Lambda}_\beta(\mathbb R^n)}\|M_\beta f\|_{M\dot{K}^{\alpha,\lambda}_{p_2,q_2}(\mathbb R^n)}
\int_{\mathbb R^n}\phi(y)|\det A^{-1}(y)|^{1/{q_2}}\sum_{j=l}^{l+m+1}2^{j(\lambda-\alpha)}dy.
\end{aligned}\end{eqnarray}

Now, for $\alpha=\lambda,$ we have

\begin{eqnarray}\begin{aligned}[b]\label{E34}\sum_{j=l}^{l+m+1}2^{j(\lambda-\alpha)}=m+2\le C(1+\log_2(\|A(y)\|\|A^{-1}(y)\|)),\end{aligned}\end{eqnarray}
and otherwise
\begin{eqnarray}\begin{aligned}[b]\label{E35}\sum_{j=l}^{l+m+1}2^{j(\lambda-\alpha)}&=2^{l(\lambda-\alpha)}\frac{1-2^{(\lambda-\alpha)(m+2)}}{1-2^{(\lambda-\alpha)}}\\
&\le C\begin{cases}\|A^{-1}(y)\|^{\alpha-\lambda}, & \alpha>\lambda,\\
\|A(y)\|^{\lambda-\alpha}, & \alpha<\lambda.\end{cases}\end{aligned}\end{eqnarray}
Hence, (\ref{E33}), (\ref{E34}) and (\ref{E35}) together yield
\begin{eqnarray*}\begin{aligned}&\|M(H_{\Phi,A}^bf)(\cdot)\|_{M\dot{K}^{\alpha,\lambda}_{p_2,q_2}(\mathbb R^n)}&\\
&\le C\|b\|_{\dot{\Lambda}_\beta(\mathbb R^n)}\|M_\beta f\|_{M\dot{K}^{\alpha,\lambda}_{p_2,q_2}(\mathbb R^n)}
\int_{\mathbb R^n}\phi(y)|\det A^{-1}(y)|^{1/{q_2}}G_{\alpha,\lambda}(y)dy.
\end{aligned}\end{eqnarray*}
Making use of the fact that $|H_{\Phi,A}^bf(x)|\le M(H_{\Phi,A}^bf)(x)$ $a.e.$ and Lemma
\ref{LMH} we get the desired result.\\

\noindent\textbf{\textit{Proof of Theorem \ref{T3}:}} Following the same arguments as given in the proof of Theorem \ref{T4}, we get

\begin{eqnarray*}\begin{aligned}&\|M(H_{\Phi,A}^bf)(\cdot)\|_{\dot{K}^{\alpha,p_2}_{q_2}(\mathbb R^n)}&\\
&\le C\|b\|_{\dot{\Lambda}_\beta(\mathbb R^n)}\|M_\beta f\|_{\dot{K}^{\alpha,p_2}_{q_2}(\mathbb R^n)}
\int_{\mathbb R^n}\phi(y)|\det A^{-1}(y)|^{1/{q_2}}\widetilde{G}_{\alpha}(y)dy.
\end{aligned}\end{eqnarray*}
However, in contrast with Theorem \ref{T4}, here we use Lemma \ref{LH} to fulfill the assertion made in the statement of this Theorem.\\

\noindent\textbf{\textit{Proof of Theorem \ref{T5}:}} For $x\in Q,$ it is easy to see that
\begin{eqnarray}\label{E37}\begin{aligned}[b]&\frac{1}{|Q|^{1+\beta/n}}\int_{Q}|H_{\Phi,A}^bf(z)-(H_{\Phi,A}^b)_Q|dz&\\
&\le\frac{2}{|Q|^{1+\beta/n}}\int_{Q}|H_{\Phi,A}^bf(z)|dz\\
&\le\frac{1}{|Q|^{1+\beta/n}}\int_{\mathbb R^n}\frac{|\Phi(y)|}{|y|^n}\int_{Q}\left|(b(z)-b_{Q})f(A(y) z)\right|dzdy\\
&+\frac{1}{|Q|^{1+\beta/n}}\int_{\mathbb R^n}\frac{|\Phi(y)|}{|y|^n}\int_{Q}\left|(b_{Q}-b_{A(y)Q})f(A(y) z)\right|dzdy\\
&+\frac{1}{|Q|^{1+\beta/n}}\int_{\mathbb R^n}\frac{|\Phi(y)|}{|y|^n}\int_{Q}\left|(b(A(y) z)-b_{A(y)Q})f(A(y) z)\right|dzdy\\
& =:J_1+J_2+J_3.\end{aligned}\end{eqnarray}

Comparing $J_1,J_2$ and $J_3$ with $I_1,I_2$ and $I_3$ estimated in the proof of Lemma \ref{LA}, one can easily estimate $J_1,J_2$ and $J_3$ by adjusting the factor
$|Q|^{\beta/n}.$ Hence, we have
\begin{eqnarray*}\begin{split}J_1\leq C\|b\|_{\dot{\Lambda}_\beta(\mathbb R^n)}\int_{\mathbb R^n}\frac{|\Phi(y)|}{|y|^n}
 M( f)(A(y)x)dy.\end{split}\end{eqnarray*}
 \begin{eqnarray*}\begin{split}J_2&\leq C\|b\|_{\dot{\Lambda}_\beta(\mathbb R^n)}\int_{\mathbb R^n}\frac{|\Phi(y)|}{|y|^n}(1+\|A(y)\|^{\beta})
 M( f)(A(y)x)dy.\end{split}\end{eqnarray*}
 \begin{eqnarray*}\begin{split}J_3&\leq C\|b\|_{\dot{\Lambda}_\beta(\mathbb R^n)}\int_{\mathbb R^n}\frac{|\Phi(y)|}{|y|^n}|\det A(y)|^{1+\beta/n}
 M( f)(A(y)x)dy.\end{split}\end{eqnarray*}

 In view of these estimates, inequality (\ref{E37}) assumes the following form
 \begin{eqnarray*}\begin{aligned}[b]&\frac{1}{|Q|^{1+\beta/n}}\int_{Q}|H_{\Phi,A}^bf(z)-(H_{\Phi,A}^b)_Q|dz&\\
&\le C\|b\|_{\dot{\Lambda}_\beta(\mathbb R^n)}\int_{\mathbb R^n}\frac{|\Phi(y)|}{|y|^n}(1+\|A(y)\|^{\beta})\max\{1,|\det A(y)|^{1+\beta/n}\}
 M( f)(A(y)x)dy.\end{aligned}\end{eqnarray*}
Applying $L^p(\mathbb R^n)$ norm on both sides and using Lemma \ref{LT}, we obtain
\begin{eqnarray*}\begin{aligned}[b]&\|H_{\Phi,A}^bf\|_{\dot{F}^{\beta,\infty}_p(\mathbb R^n)}\le C\|b\|_{\dot{\Lambda}_\beta(\mathbb R^n)}&\\
&\qquad\times\int_{\mathbb R^n}\frac{|\Phi(y)|}{|y|^n}(1+\|A(y)\|^{\beta})\max\{1,|\det A(y)|^{1+\beta/n}\}
 \|Mf(A(y)\cdot)\|_{L^p(\mathbb R^n)}dy.\end{aligned}\end{eqnarray*}

Finally by scaling argument and boundedness of $M$ on $L^p(\mathbb R^n),$ (see \cite{GF}) we have
\begin{eqnarray*}\begin{split}\|H_{\Phi,A}^bf\|_{F^{\beta,\infty}_p(\mathbb R^n)}\le C K_5\|b\|_{\dot{\Lambda}_\beta(\mathbb R^n)}\|f\|_{L^{p}({\mathbb R^n})},\end{split}\end{eqnarray*}
which is as required.

\section{\textbf{Central-$BMO$ Estimates for $H_{\Phi,A}^b$ on Herz-type Space}}

\subsection{Main Results}
\begin{theorem}\label{T6}
Let $1<p,q_1,q_2<\infty,1/q_2=1/q+1/q_1,\lambda>0$ and $\alpha_2\in \mathbb R.$  If $\alpha_1=n/q+\alpha_2$ and $b\in C\dot{M}O^q(\mathbb R^n),$ then
$$\|H_{\Phi,A}^bf\|_{M\dot{K}^{\alpha_2,\lambda}_{p,q_2}(\mathbb R^n)}\le CK_6\|b\|_{C\dot{M}O^q(\mathbb R^n)}\|f\|_{M\dot{K}^{\alpha_1,\lambda}_{p,q_1}(\mathbb R^n)},$$
where
\begin{eqnarray*}\begin{aligned}[b] K_6&=\int_{\mathbb R^n}\frac{|\Phi(y)|}{|y|^n}|\det A^{-1}(y)|^{1/q_1}G_{\alpha_1,\lambda}(y)\\
&\quad\times\left(\log\frac{2}{\|A(y)\|}\chi_{\{\|A(y)\|<1\}}+\log2\|A(y)\|\chi_{\{\|A(y)\|\ge1\}}\right)dy,\end{aligned}\end{eqnarray*}
and $G_{\alpha_1,\lambda}(y)$ is the same function as given in (\ref{GLE}) with $\alpha$ is replaced by $\alpha_1.$

\end{theorem}
\begin{theorem}\label{T7}
Let $1<p,q_1,q_2<\infty,1/q_2=1/q+1/q_1$ and $\alpha_2\in \mathbb R.$  If $\alpha_1=n/q+\alpha_2$ and $b\in C\dot{M}O^q(\mathbb R^n),$ then
$$\|H_{\Phi,A}^bf\|_{\dot{K}^{\alpha_2,p}_{q_2}(\mathbb R^n)}\le CK_7\|b\|_{C\dot{M}O^q(\mathbb R^n)}\|f\|_{\dot{K}^{\alpha_1,p}_{q_1}(\mathbb R^n)},$$
where
\begin{eqnarray*}\begin{aligned}[b] K_7&=\int_{\mathbb R^n}\frac{|\Phi(y)|}{|y|^n}|\det A^{-1}(y)|^{1/q_1}\widetilde{G}_{\alpha_1}(y)\\
&\quad\times\left(\log\frac{2}{\|A(y)\|}\chi_{\{\|A(y)\|<1\}}+\log2\|A(y)\|\chi_{\{\|A(y)\|\ge1\}}\right)dy,\end{aligned}\end{eqnarray*}
and $\widetilde{G}_{\alpha_1}(y)$ is the same function as given in (\ref{GE}) with $\alpha$ is replaced by $\alpha_1.$

\end{theorem}

\subsection{\textbf{Proof of Theorem \ref{T7}.}} Here, we decompose $\|H_{\Phi,A}^bf\|_{L^{q_2}(C_k)}$ as:
\begin{eqnarray*}\begin{aligned}[b]
\|H_{\Phi,A}^bf\|_{L^{q_2}(C_k)}&=\left(\int_{C_k}\left|\int_{\mathbb R^n}\frac{\Phi(y)}{|y|^n}(b(x)-b(A(y) x))f(A(y) x)dy\right|^{q_2}dx\right)^{1/q_2}\\
&\le\int_{\mathbb R^n}\frac{|\Phi(y)|}{|y|^n}\left(\int_{C_k}\left|(b(x)-b(A(y) x))f(A(y) x)\right|^{q_2}dx\right)^{1/q_2}dy\\
&\le\int_{\mathbb R^n}\frac{|\Phi(y)|}{|y|^n}\left(\int_{C_k}\left|(b(x)-b_{B_k})f(A(y) x)\right|^{q_2}dx\right)^{1/q_2}dy\\
&+\int_{\mathbb R^n}\frac{|\Phi(y)|}{|y|^n}\left(\int_{C_k}\left|(b_{B_k}-b_{A(y)B_k})f(A(y) x)\right|^{q_2}dx\right)^{1/q_2}dy\\
&+\int_{\mathbb R^n}\frac{|\Phi(y)|}{|y|^n}\left(\int_{C_k}\left|(b(A(y) x)-b_{A(y)B_k})f(A(y) x)\right|^{q_2}dx\right)^{1/q_2}dy\\&= L_1+L_2+L_3,
\end{aligned}\end{eqnarray*}
By H\"older inequality and change of variables it is simple to have
\begin{eqnarray*}
    \begin{aligned}[b]
L_1&\le\int_{\mathbb R^n}\frac{|\Phi(y)|}{|y|^n}\left(\int_{C_k}\left|b(x)-b_{B_k}\right|^{q} dx\right)^{1/q}
\left(\int_{C_k}\left|f(A(y) x)\right|^{q_1} dx\right)^{1/q_1}dy\\
&\le |B_k|^{1/q}\|b\|_{C\dot{M}O^{q}(\mathbb R^n)}\int_{\mathbb R^n}\frac{|\Phi(y)|}{|y|^n}|\det A^{-1}(y)|^{1/q_1}\|f\|_{L^{q_1}(A(y)C_k)}dy.
\end{aligned}\end{eqnarray*}

In order to estimate $L_2,$ we rewrite it as
\begin{eqnarray}\label{L2}\begin{aligned}[b]
L_2&=\int_{\mathbb R^n}\frac{|\Phi(y)|}{|y|^n}|\det A^{-1}(y)|^{1/q_2}\|f\|_{L^{q_2}(A(y)C_k)}\left|b_{B_k}-b_{A(y)B_k}\right|dy\\
&\le\int_{\mathbb R^n}\frac{|\Phi(y)|}{|y|^n}|\det A^{-1}(y)|^{1/q_2}|A(y)B_k|^{1/q}\|f\|_{L^{q_1}(A(y)C_k)}\left|b_{B_k}-b_{A(y)B_k}\right|dy\\
&=|B_k|^{1/q}\int_{\mathbb R^n}\frac{|\Phi(y)|}{|y|^n}|\det A^{-1}(y)|^{1/q_1}\|f\|_{L^{q_1}(A(y)C_k)}\left|b_{B_k}-b_{A(y)B_k}\right|dy\\
& =|B_k|^{1/q} \int_{\|A(y)\|<1}\frac{|\Phi(y)|}{|y|^n}|\det A^{-1}(y)|^{1/q_1}\|f\|_{L^{q_1}(A(y)C_k)}\left|b_{B_k}-b_{A(y)B_k}\right|dy \\
               &\quad+|B_k|^{1/q}\int_{\|A(y)\|\ge1}\frac{|\Phi(y)|}{|y|^n}|\det A^{-1}(y)|^{1/q_1}\|f\|_{L^{q_1}(A(y)C_k)}
               \left|b_{B_k}-b_{A(y)B_k}\right|dy \\
               &=: |B_k|^{1/q}\left( L_{21} +L_{22} \right).
\end{aligned}\end{eqnarray}

Thus, for $\|A(y)\|<1,$
we have
\begin{equation*}
    \begin{split}
    &L_{21}=\sum_{j=0}^\infty\int_{2^{-j-1}\le\|A(y)\|<2^{-j}}\frac{|\Phi(y)|}{|y|^n}|\det A^{-1}(y)|^{1/q_1}\|f\|_{L^{q_1}(A(y)C_k)}\\
    &\qquad\qquad\qquad\times\left\{\sum_{i=1}^j|b_{2^{-i}B_k}-b_{2^{-i+1}B_k}|
    +|b_{2^{-j}B_k}-b_{A(y)B_k}|\right\}dy.\end{split}
\end{equation*}
A use of H\"older inequality yields
\begin{equation*}
    \begin{split}
    |b_{2^{-i}B_k}-b_{2^{-i+1}B_k}|&\le \frac{1}{|2^{-i}B_k|}\int_{2^{-i}B_k}|b(y)-b_{2^{-i+1}B_k}|dy\\
    &\le \frac{|2^{-i+1}B_k|}{|2^{-i}B_k|}\|b\|_{C\dot{M}O^{q}(\mathbb R^n)}\\
    &\le C \|b\|_{C\dot{M}O^{q}(\mathbb R^n)}.\end{split}
\end{equation*}
Similarly, $|b_{2^{-j}B_k}-b_{A(y)B_k}|\le C \|b\|_{C\dot{M}O^{q}(\mathbb R^n)}$ and therefore
\begin{equation*}
    \begin{split}
    L_{21}&\le C \|b\|_{C\dot{M}O^{q}(\mathbb R^n)}\\
    &\quad\times\sum_{j=0}^\infty\int_{2^{-j-1}\le\|A(y)\|<2^{-j}}\frac{|\Phi(y)|}{|y|^n}|\det A^{-1}(y)|^{1/q_1}\|f\|_{L^{q_1}(A(y)C_k)}\left\{j+1\right\}dy\\
    &\le C \|b\|_{C\dot{M}O^{q}(\mathbb R^n)}\\
    &\quad\times\sum_{j=0}^\infty\int_{2^{-j-1}\le\|A(y)\|<2^{-j}}\frac{|\Phi(y)|}{|y|^n}|\det A^{-1}(y)|^{1/q_1}\|f\|_{L^{q_1}(A(y)C_k)}\left\{\log2^j+1\right\}dy\\
    &\le C \|b\|_{C\dot{M}O^{q}(\mathbb R^n)}\int_{\|A(y)\|<1}\frac{|\Phi(y)|}{|y|^n}|\det A^{-1}(y)|^{1/q_1}\|f\|_{L^{q_1}(A(y)C_k)}\log\frac{2}{\|A(y\|}dy.
    \end{split}
\end{equation*}

Similar Arguments result in the following estimation of $L_{22}.$
\begin{equation*}
    \begin{split}
    L_{22}&\le C \|b\|_{C\dot{M}O^{q}(\mathbb R^n)}\int_{\|A(y)\|\ge1}\frac{|\Phi(y)|}{|y|^n}|\det A^{-1}(y)|^{1/q_1}\|f\|_{L^{q_1}(A(y)C_k)}\log2\|A(y\|dy.
    \end{split}
\end{equation*}
Having these estimates of $L_{21}$ and $L_{22},$ (\ref{L2}) assumes the following form
\begin{equation*}
    \begin{aligned}[b]
    L_2\le C  & |B_k|^{1/q}\|b\|_{C\dot{M}O^{q}(\mathbb R^n)}\\
    &{}  \Bigg( \int_{\|A(y)\|<1}\frac{|\Phi(y)|}{|y|^n}|\det A^{-1}(y)|^{1/q_1}\|f\|_{L^{q_1}(A(y)C_k)}\log\frac{2}{\|A(y\|}dy \\
               & ~ +\int_{\|A(y)\|\ge1}\frac{|\Phi(y)|}{|y|^n}|\det A^{-1}(y)|^{1/q_1}\|f\|_{L^{q_1}(A(y)C_k)}\log2\|A(y\|dy \Bigg).
    \end{aligned}
\end{equation*}

It remains to estimate $L_3.$ For this purpose we take advantage of H\"older inequality to obtain
\begin{eqnarray*}\label{CEE}\begin{split}
&\|(b(A(y)\cdot)-b_{A(y)B_k})f(A(y)\cdot )\|_{L^{q_2}(C_k)}\\&=\left(\int_{C_k}\left|(b(A(y)x)-b_{A(y)B_k})f(A(y) x)\right|^{q_2} dx\right)^{1/q_2}\\
&=|\det A^{-1}(y)|^{1/q_2}\left(\int_{A(y)C_k}\left|(b(x)-b_{A(y)B_k})f( x)\right|^{q_2} dx\right)^{1/q_2}\\
&\le |\det A^{-1}(y)|^{1/q_2}
\|b-b_{A(y)B_k}\|_{L^{q}(A(y)C_k)}\|f\|_{L^{q_1}(A(y)C_k)}\\
&\le    |\det A^{-1}(y)|^{1/q_2}|A(y)B_k|^{1/q}
\|b\|_{C\dot{M}O^{q}(\mathbb R^n)}\|f\|_{L^{q_1}(A(y)C_k)}\\
&=|B_k|^{1/q}|\det A^{-1}(y)|^{1/q_1}
\|b\|_{C\dot{M}O^{q}(\mathbb R^n)}\|f\|_{L^{q_1}(A(y)C_k)}.
\end{split}\end{eqnarray*}
Hence,
\begin{eqnarray*}\begin{split}
L_3&\le |B_k|^{1/q}\|b\|_{C\dot{M}O^{q}(\mathbb R^n)}\int_{\mathbb R^n}\frac{|\Phi(y)|}{|y|^n}|\det A^{-1}(y)|^{1/q_1}\|f\|_{L^{q_1}(A(y)C_k)}dy.
\end{split}\end{eqnarray*}

Combining the estimates for $L_1,$ $L_2,$ and $L_3$ we obtain
\begin{eqnarray*}\begin{split}\|H_{\Phi,A}^bf\|_{L^{q_2}(C_k)}&\le C|B_k|^{1/q}\|b\|_{C\dot{M}O^{q}(\mathbb R^n)}\int_{\mathbb R^n}\frac{|\Phi(y)|}{|y|^n}|\det A^{-1}(y)|^{1/q_1}\|f\|_{L^{q_1}(A(y)C_k)}\\
&\qquad\times\left(\log\frac{2}{\|A(y)\|}\chi_{\{\|A(y)\|<1\}}+\log2\|A(y)\|\chi_{\{\|A(y)\|\ge1\}}\right)dy.\end{split}\end{eqnarray*}
 Again, to make our calculation convenient, we use the following notation $$\varphi(y)=\frac{|\Phi(y)|}{|y|^n}|\det A^{-1}(y)|^{1/q_1}
\left(\log\frac{2}{\|A(y)\|}\chi_{\{\|A(y)\|<1\}}+\log2\|A(y)\|\chi_{\{\|A(y)\|\ge1\}}\right).$$ Thus, we rewrite above inequality as
\begin{eqnarray}\label{SE}\begin{split}\|H_{\Phi,A}^bf\|_{L^{q_2}(C_k)}&\le C|B_k|^{1/q}
\|b\|_{C\dot{M}O^{q}(\mathbb R^n)}\int_{\mathbb R^n}\varphi(y)\|f\|_{L^{q_1}(A(y)C_k)}dy.\end{split}\end{eqnarray}

 Still we have to approximate $\|f\|_{L^{q_1}(A(y)C_k)}.$ For this end we infer from (\ref{EA}) that
 \begin{eqnarray*}\begin{aligned}
\|f\|_{L^{q_1}(A(y)C_k)}\le\sum_{j=l}^{l+m+1}\|f\|_{L^{q_1}(C_{k+j})}.
\end{aligned}
\end{eqnarray*}
By this inequality, (\ref{SE}) becomes
\begin{eqnarray}\label{SE1}\begin{split}\|H_{\Phi,A}^bf\|_{L^{q_2}(C_k)}&\le C|B_k|^{1/q}
\|b\|_{C\dot{M}O^{q}(\mathbb R^n)}\int_{\mathbb R^n}\varphi(y)\sum_{j=l}^{l+m+1}\|f\|_{L^{q_1}(C_{k+j})}dy.\end{split}\end{eqnarray}

 Next by definition of Morrey-Herz space $M\dot{K}^{\alpha_2,\lambda}_{p,q_2}(\mathbb R^n)$, (\ref{SE1}), Minkowski inequality and the condition $\alpha_1=\alpha_2+n/q,$ we get
\begin{eqnarray*}\begin{aligned}[b] &\|H_{\Phi,A}^bf\|_{M\dot{K}^{\alpha_2,\lambda}_{p,q_2}(\mathbb R^n)}&\\
&=\sup_{k_0\in\mathbb Z}2^{-k_0\lambda}
\left\{\sum_{k=-\infty}^{k_0} 2^{k\alpha_2 p}\|H_{\Phi,A}^bf\|_{L^{q_2}(C_k)}^p\right\}^{1/p}\\
&\le C\|b\|_{C\dot{M}O^q(\mathbb R^n)}
\int_{\mathbb R^n}\varphi(y)\\
&\quad\times\sup_{k_{0}\in \mathbb Z}2^{-k_0\lambda}\left\{\sum^{k_{0}}_{k=-\infty} \left(\sum_{j=l}^{l+m+1}2^{k(\alpha_2+n/q)}\|f\|_{L^{q_1}(C_{k+j})}\right)^p\right\}^{1/p}dy\\
&\le C\|b\|_{C\dot{M}O^q(\mathbb R^n)}
\int_{\mathbb R^n}\varphi(y)\\
&\quad\times\sup_{k_{0}\in \mathbb Z}2^{-k_0\lambda}
\sum_{j=l}^{l+m+1}2^{-j\alpha_1}\left\{\sum^{k_{0}+j}_{k=-\infty} 2^{k\alpha_1 p }\| f\|_{L^{q_1}(C_{k})}^p\right\}^{1/p}dy\\
&\le C\|b\|_{C\dot{M}O^q(\mathbb R^n)}\|f\|_{M\dot{K}^{\alpha_1,\lambda}_{p,q_2}(\mathbb R^n)}
\int_{\mathbb R^n}\varphi(y)\sum_{j=l}^{l+m+1}2^{j(\lambda-\alpha_1)}dy\\
&\le C\|b\|_{C\dot{M}O^q(\mathbb R^n)}\|f\|_{M\dot{K}^{\alpha_1,\lambda}_{p,q_2}(\mathbb R^n)}
\int_{\mathbb R^n}\varphi(y)G_{\alpha_1,\lambda}(y)dy,
\end{aligned}\end{eqnarray*}
where $G_{\alpha_1,\lambda}(y)$ is the same function as given in (\ref{GLE}) with $\alpha$ is replaced by $\alpha_1.$ This finishes the proof of Theorem \ref{T6}.\\

\subsection{\textbf{Proof of Theorem \ref{T7}.}}
Since, the the proof of Theorem \ref{T6} and \ref{T7} are symmetrical. So, by definition of Herz space $\dot{K}^{\alpha_2,p}_{q_2}(\mathbb R^n)$, (\ref{SE1}), Minkowski inequality and the condition $\alpha_1=\alpha_2+n/q,$ we get
\begin{eqnarray*}\begin{aligned}[b] &\|H_{\Phi,A}^bf\|_{\dot{K}^{\alpha_2,p}_{q_2}(\mathbb R^n)}
&\le C\|b\|_{C\dot{M}O^q(\mathbb R^n)}\|f\|_{\dot{K}^{\alpha_1,p}_{q_1}(\mathbb R^n)}
\int_{\mathbb R^n}\varphi(y)\widetilde{G}_{\alpha_1}(y)dy,
\end{aligned}\end{eqnarray*}
where $\widetilde{G}_{\alpha_1}(y)$ is the same function as given in (\ref{GE}) with $\alpha$ is replaced by $\alpha_1.$ This completes the proof of Theorem \ref{T7}.\\


\begin{thebibliography}{99}

\bibitem{BS}    A. Baernstein and E. Sawyer, Embedding and multiplier theorems for $H^p(\mathbb R^n),$ Mem. Amer. Math. Soc. {\bf{53}} (1985), 1--82.

\bibitem{BGG}   V. Burenkov, H. Guliyev and V. Guliyev, Necessary and sufficient conditions for the bounndedness of fractional maximal operators
                in local Morrey type spaces, J. Comput. Appl. Math. {\bf{208}} (2007), 280--301.

\bibitem{BL}    V. Burenkov and E. Liflyand, On the boundedness of Hausdorff operators on Morrey-type spaces, Eurasian Math. J. {\bf{8}} (2017), 97--104.

\bibitem{CFL}   J. Chen, D. Fan and J. Li,  Hausdorff operators on function spaces, Chin. Ann. Math. {\bf 33} (2012), 537--556.

\bibitem{CHZ}   J. Chen, S. He and X. Zhu, Boundedness of Hausdorff operators on the power weighted Hardy spaces, Appl. Math. {\bf{32}} (2017) 462--476.

\bibitem{CFW}   J. Chen, D. Fan and S. Wang, Hausdorff operators on Eucleadian spaces, Appl. Math. {\bf{28}} (2013), 548--564.

\bibitem{CZ}   J. Chen and X. Zhu, Boundedness of multidimensional Hausdorff operators on $H^1(\mathbb R^n),$
                 J. Math. Anal. Appl. {\bf{409}} (2014), 428--434.
\bibitem{CG}   M. Christ and L. Grafakos, Best constants for two nonconvolution inequalities, Proc. Amer. Math. Soc. {\bf{123}} (1995), 1687-1693.

\bibitem{GJ}    G. Gao and H. Jia, Boundedness of commutators of high dimensional Hausdorff operator.
                J. Function Spaces Appl. {\bf{2012}} (2012):54120.

\bibitem{GWAG} G. Gao, X. Wu, A. Hussain and G. Zhao, Some esimates for Hausdorff operators, J. Math. Inequal. {\bf{9}} (2015), 641--651.

\bibitem{GZ}   G. Gao and Y. Zhong, Some inequalities for Hausdorff operators, Math. Ineq. Appl. {\bf{17}} (2014), 1061--1078.

\bibitem{GF}    L. Grafakos, Classical Fourier Analysis, second edition, Graduate Texts in Mathematics, 249, Springer, New York, 2008.

\bibitem{GU}    Y. Guo, Boundedness of some operators on non-homogeneous spaces, Master Degree Dissertation, Beijing Normal University.

\bibitem{GZ1}    J. Guo and F. Zhao, Some $q$-inequalities for Hausdorff operators, Front. Math. China, {\bf 12} (2017),  879--889.

\bibitem{KP}    K. Ho, Hardy Littlewood P\'olya inequalities and Hausdorff operator on Block spaces, Math. Inequal. Appl., {\bf{19}} (2016), 697--707.

\bibitem{HKQ}   H. Hung, L. Ky and T. Quang, Norm of the Hausdorff Operator on the Real Hardy Space $H^1(\mathbb R),$
                Complex Anal. Oper. Theory {\bf{12}} (2018), 235--245.

\bibitem{HA}   A. Hussain and A. Ajaib, Some weighted inequalities for Hausdorff operators and commutators, J. Inequal. Appl.  {\bf{2018}} (2018)6, 19 pages.

\bibitem{HA1}    A. Hussain and M. Ahmed, Weak and strong type estimates for the commutators of Hausdorff operator, Math. Inequal. Appl. {\bf{20}} (2017), 49--56.

\bibitem{HG1}   A. Hussain and G. Gao, Some new estimates for the commutators of n-dimensional Hausdorff operator,
                Appl. Math. {\bf{29}} (2014), 139--150.

\bibitem{HG2}   A. Hussain and G. Gao, Multidimensional Hausdorff operators and commutators on Herz-type spaces,
                J. Ineq. Appl. {\bf{2013}} (2013)594, 12 pages.

\bibitem{K}     Y. Kanjin, The Hausdorff operators on the real Hardy spaces $H^p(\mathbb R),$ Studia Math. {\bf 148} (2001), 37--45.

\bibitem{LL}    A. Lerner and E. Liflyand, Multidimensional Hausdorff operators on the real Hardy spaces, Jour. Austr. Math. Soc.
                {\bf 83} (2007), 79--86.

\bibitem{EL2}   E. Lilyand, Boundedness of multidimensional Hausdors operators on $H^1(\mathbb R^n)$, Acta Sci. Math. (Szeged) {\bf{74}} (2008), 845--851.

\bibitem{EL1}   E. Liflyand, Hausdorff operators on Hardy spaces, Eurasian Math. J. {\bf{4}} (2013), 101--141.

\bibitem{LM}    E. Liflyand and A. Miyachi, Boundedness of the Hausdorff operators in $H^{p}$ spaces, $0<p<1,$
                Stud. Math. {\bf{194}} (2009), 279--292.

\bibitem{LM1}   E. Liflyand and F. M\'{o}recz, The Hausdorff operator is bounded on the real Hardy space $H^{1}(\mathbb R)$, Proc. Amer. Math. Soc. {\bf 128} (2000), 1391--1396.

\bibitem{LY}    S. Lu and D. Yang, Herz-type Sobolev and Bessel potential spaces and their applications, Sci. China Ser. A {\bf{40}} (1997), 113--129.

\bibitem{LYH}   S. Lu, D. Yang and G. Hu, Herz Type Spaces and Their Application, Science Press-Beijing, China, 2008.

\bibitem{CM}   C. Morrey, On the solutions of quasi-linear elliptic partial differential equations,
                Trans. Amer. Math. Soc. {\bf{43}} (1938), 126--166.

\bibitem{PM}   M. Paluszynski, Characterization of Besov spaces via the commutator operator of Coifman, Rochberg and Weiss,
                Ind. Univ. Math. J. {\bf{44}} (1995), 1--18.

\bibitem{RF2}   J. Ruan and D. Fan, Hausdorff operators on the power weighted Hardy spaces, J. Math. Anal. Appl. {\bf{433}} (2016), 31--48.

\bibitem{RF1}   J. Ruan and D. Fan, Hausdorff operators on weighted Herz-type Hardy spaces, Math. Inequal. Appl. {\bf{19}} (2016), 565–-587.

\bibitem{RFW}   J. Ruan, D. Fan, and Q. Wu, Weighted Herz space estimates for Hausdorff operators on the Heisenberg group,
                Banach J. Math. Anal.  {\bf{11}} (2017), 513--535.

\bibitem{XW}    X. Wu, Necessary and sufficient conditions for generalized Hausdorff operators and commutators, Ann. Funct. Anal. {\bf{6}} (2015) 60--72.







\end{thebibliography}
\end{document}